\documentclass[11pt]{amsart}

\usepackage{amsthm, xcolor}
\usepackage{latexsym, multicol}
\usepackage{thmtools, thm-restate}

\usepackage{ulem}

\title{Skew Zero Forcing Matroids}

\thanks{Research was sponsored by the Air Force Office of Scientific Research and was accomplished under Grant Number W911NF-25-1-0110. The views and conclusions contained in this document are those of the authors and should not be interpreted as representing the official policies, either expressed or implied, of the Army Research Office or the U.S. Government. The U.S. Government is authorized to reproduce and distribute reprints for Government purposes notwithstanding any copyright notation herein.}

\author{Dennis Hall}
\address{Department of Mathematics, Angelo State University, San Angelo, Texas, USA}
\email{Dennis.Hall@angelo.edu}

\author{Simon Pfeil}
\address{Department of Mathematics, Angelo State University, San Angelo, Texas, USA}
\email{Simon.Pfeil@angelo.edu}

\author{Jesse Taylor}
\address{Department of Mathematics, Angelo State University, San Angelo, Texas, USA}
\email{Jesse.Taylor@angelo.edu}

\usepackage{amsmath}
\usepackage{amssymb}

\usepackage{graphicx}
\usepackage{subcaption}
\usepackage{comment}

\usepackage{url}
\usepackage[margin=1in]{geometry}
\usepackage[table]{xcolor}
\usepackage{cancel} % allows strikethrough in math mode

\definecolor{jesse}{HTML}{128A2C}

\definecolor{dennis}{HTML}{7030A0}

\definecolor{simon}{HTML}{1679c9}

\usepackage[skins,breakable]{tcolorbox}

\newtcolorbox{dennis-stblock}{
  blanker,
  breakable,
  parbox=false,
  finish={
    \draw[dennis,line width=0.8pt]
      (frame.south west) -- (frame.north east);
    \draw[dennis,line width=0.8pt]
      (frame.north west) -- (frame.south east);
  }
}

\newtcolorbox{simon-stblock}{
  blanker,
  breakable,
  parbox=false,
  finish={
    \draw[simon,line width=0.8pt]
      (frame.south west) -- (frame.north east);
    \draw[simon,line width=0.8pt]
      (frame.north west) -- (frame.south east);
  }
}

\usepackage{pgfplots}
\pgfplotsset{compat=1.18}
\usepackage{tikz}
\usetikzlibrary{calc, intersections, automata, shapes.geometric}

\newtheorem{theorem}{Theorem}[section]
\newtheorem{lemma}[theorem]{Lemma}
\newtheorem{corollary}[theorem]{Corollary}
\newtheorem{proposition}[theorem]{Proposition}

\theoremstyle{definition}
\newtheorem{definition}[theorem]{Definition}

\theoremstyle{remark}

\newcommand{\myloop}[2]{\draw (#1) to[out=#2, in=#2+60, looseness=25] (#1)}

\newcommand{\cl}{\mathop{\mathrm{cl}}}

\numberwithin{equation}{section}

\newcommand{\ba}{\backslash}

\begin{document}

\begin{abstract}
Matroid theory and skew zero forcing (SZF) in graphs are both topics of
mathematical interest, but very little is known about their overlap. When a matroid arises from a graph in the traditional way, the edges of the graph correspond with the elements of the matroid. However, the SZF color change rule for graphs can be viewed as a closure operator on the vertices of a graph. As such, when a matroid arises from a graph via SZF, the vertices of the graph correspond with the elements of the matroid. We call such a matroid SZF-realizable. Notably, unlike the traditional method of obtaining a matroid from a graph, skew zero forcing in a graph does not always give rise to a matroid. This paper establishes fundamental results in graph and matroid theory with respect to (1) which matroids are SZF-realizable, and (2) what kinds of graphs do and do not give rise to SZF matroids. The main result of the paper focuses on cycle matroids of outerplanar graphs.
\end{abstract}

\maketitle

\section{Introduction}
\label{S.intro}

This paper explores the intersection of matroids and the skew zero forcing (SZF) color change rule in graphs. SZF in graphs was first considered in \cite{allison}, and it has since received significant interest from the mathematics community. In March 2023, a preprint \cite{cooper} looked at SZF as a closure operator and, more specifically, observed that SZF sometimes, but not always, gives rise to a matroid closure operator. Additionally, it established a polynomial time algorithm for calculating the skew zero forcing number of a graph when the graph's SZF closure operator is also a matroid closure operator.

We assume the reader is familiar with fundamental definitions of matroid theory. Except where otherwise noted, the notation and definitions in this paper follow \cite{oxley}. The way that matroids arise from SZF as a closure operator is unlike the typical way that matroids arise from graphs. In particular, the cycle matroid that arises from a graph uses the graph's edge set as the ground set of the matroid. However, when using an SZF closure operator, the vertex set of the graph corresponds to the ground set of the matroid. 

Let $G = (V,E)$ be a graph. The vertices of $G$ will each be labeled ``filled'' or ``unfilled," and a vertex is unfilled unless otherwise specified. A subset $S \subseteq V$ defines an \textit{initial coloring} by filling all the vertices of $S$ and leaving all the vertices not in $S$ unfilled. Unfilled vertices in $G$ can be filled after the initial coloring via the \textit{skew color change rule} as described in \cite{allison}, which works in the following way. For every vertex $v\in V$, consider its neighborhood, $n(v)$. If there is exactly one unfilled vertex in $n(v)$, then fill that vertex; otherwise, do nothing. Repeat this until there are no neighborhoods with exactly one unfilled vertex. The resulting set of filled vertices is considered closed; that is, a set of vertices that does not change the filled status of any vertex in $V$ after applying the skew color change rule. We say vertices in the initial set $S$ and vertices filled in by the skew color change rule are \textit{forced by S}. We say a vertex $v$ \textit{forces} a vertex $u$ in $G$, if the skew color change rule results in $u$ being filled because it was the only unfilled neighbor of $v$. The \textit{skew zero forcing closure} of an initial filled set $S$ is the result of applying the skew color change rule until no more change occurs. We denote this closure by $\textrm{scl}_G(S)$, or simply $\textrm{scl}(S)$, if there is no confusion surrounding the graph.

Unlike most researchers considering skew zero forcing in graphs (\cite{AIM},\cite{fail},\cite{rank}, for example), we do not insist that our graphs are simple. In our expanded SZF definition, we work on multigraphs and apply the skew color change rule using adjacency {\it counted with multiplicity}. A vertex $v$ (filled or not) may force $u$ precisely when $u$ is the only unfilled neighbor of $v$ after counting incidences. A loop at $v$, for example, counts $v$ as its own neighbor; thus a single loop lets an unfilled vertex $v$ self-force once all other neighbors are filled. A double loop on an unfilled vertex blocks self-forcing since a double-looped vertex will count itself as two unfilled neighbors. In fact, a double-looped vertex acts as a gate for forcing, allowing forcing through that vertex only once it is filled itself. Parallel edges play a similar role between distinct vertices: multiple edges from $u$ to $v$ means $v$ would count $u$ as multiple neighbors when applying the skew color change rule. One reason for this expansion to multigraphs is addressed in Section \ref{S.minors}, where the minor operation for both graphs and matroids is explored with respect to SZF.

From \cite{oxley}, a matroid closure operator $\textrm{cl}:\mathcal{P}(E)\rightarrow \mathcal{P}(E)$ on a ground set $E$ has the following properties.
\begin{enumerate}
\item[(CL1)] If $X \subseteq E$, then $X \subseteq \cl(X)$.
\item[(CL2)] If $X \subseteq Y \subseteq E$, then $\cl(X) \subseteq \cl(Y)$.
\item[(CL3)] If $X \subseteq E$, then $\cl(\cl(X)) = \cl(X)$.
\item[(CL4)] If $X \subseteq E$ and $x\in E$, and $y \in \cl(X\cup x) - \cl(X)$, then $x\in \cl(X\cup y)$.
\end{enumerate}

Similar to above, we denote this closure $\cl_M(X)$, or simply $\cl(X)$, when there is no confusion regarding the matroid. An SZF closure operator of a graph always satisfies (CL1)-(CL3). When it also satisfies (CL4), then the SZF closure operator is a matroid closure operator. In this case, we say the graph is \textit{SZF-matroidal}, and denote its corresponding matroid $SZF(G)$. Similarly, when a matroid $M$ can be constructed via SZF-closure on some graph, we call $M$ \textit{SZF-realizable}; however, when a matroid $M$ is SZF-realizable there may be several non-isomorphic graphs that give rise to $M$. Thus, when a matroid $M$ is SZF-realizable we denote the set of all graphs that correspond to $M$ by $M_{SZF}$. When it is sufficient to refer to any one member of $M_{SZF}$ and the context is clear, we use $M_{SZF}$ to mean some single member of that set.

Graphs exist that are not SZF-matroidal and matroids exist that are not SZF-realizable. These two facts combine to establish a relationship between SZF and matroids whose details are unclear. This paper seeks clarity on this relationship. In Section \ref{S.main}, we focus on the specific case of graphic matroids; that is, matroids that can be realized as the cycle matroid of some graph.

In \cite{cooper} it is shown that for any set $S$, the closure $\textrm{scl}_G(S)$ is unique and respects subset containment. Note that if a graph $G$ is SZF-matroidal, then the matroid $SZF(G)$ is unique. However, as mentioned above, if a matroid $M$ is SZF-realizable, there may be many graphs on which the SZF closure operator gives rise to $M$. Section \ref{S.prelims} presents fundamental results for handling the interaction between graphs and matroids via SZF, focusing primarily on how graph structure and the corresponding matroid structure are related. 

\begin{figure}[ht]
\begin{center}
\begin{tikzpicture}[scale=1.4,
  vertex/.style={circle,draw,fill=white,inner sep=1.5pt}
]

% P_3
\node[vertex] (v1) at  (1,1) {};
\node[vertex] (v2) at  (2,1) {};
\node[vertex] (v3) at  (3,1) {};

\draw (v1) -- (v2) -- (v3);

% U_2,3
\draw (5,1) -- (6,1) -- (7,1);
\draw[fill] (5,1) circle (0.1);
\draw[fill] (6,1) circle (0.1);
\draw[fill] (7,1) circle (0.1);

\end{tikzpicture}
\end{center}

\caption{The graph consisting of a path with three vertices (left) and the geometric representation of the matroid $U_{2,3}$ (right).}
\label{F.versus}
\end{figure}
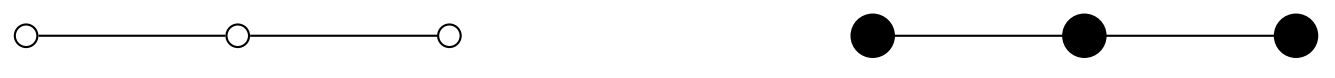

The class of uniform matroids will be useful in the discussion that follows. We use $U_{r,n}$,  where $r$ and $n$ are integers and $0\leq r\leq n$, to denote the uniform matroid of rank $r$ on $n$ elements. Throughout this paper, there are pictures of graphs and geometric representations of matroids. To avoid confusion, graphs will be drawn with white (unfilled) vertices, while geometric representations of matroids will be drawn with elements as filled-in circles, as in Figure \ref{F.versus}.

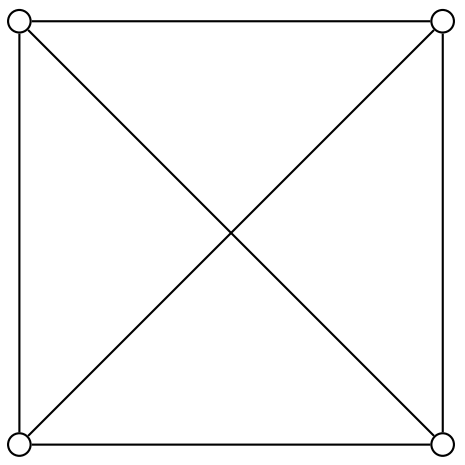
\begin{figure}[ht]
\begin{center}
\begin{tikzpicture}[scale=1.4,
  vertex/.style={circle,draw,fill=white,inner sep=1.5pt}
]

%K4 non-embedding
\node[vertex] (v1) at  (2,-3) {};
\node[vertex] (v2) at  (4,-3) {};
\node[vertex] (v3) at  (2,-1) {};
\node[vertex] (v4) at  (4,-1) {};

\draw (v1) -- (v2) -- (v4) -- (v3) -- (v1);
\draw (v1) -- (v4);
\draw (v2) -- (v3);

\end{tikzpicture}
\end{center}

\caption{The complete graph on four vertices, $K_4$.}
\label{k4 graph}
\end{figure}

We conclude this section with an example of a graph that gives rise to a uniform matroid via the SZF closure operator, and compare it to the cycle matroid of that graph. Consider the complete graph on four vertices, $K_4$, depicted in Figure \ref{k4 graph}. Observe that filling any one vertex will not force any unfilled vertices, but that filling any two vertices will force the entire graph. Thus, $SZF(K_4) \cong U_{2,4}$, whose geometric representation is shown in Figure \ref{u24}. We note that the typical way to obtain a matroid from $K_4$ would be its cycle matroid $M(K_4)$, whose geometric representation is shown in Figure \ref{k4 matroid}.

\begin{figure}[ht]
\begin{center}
\begin{tikzpicture}

%U_2,4
\draw (1,1) -- (2,1) -- (3,1) -- (4,1);
\draw[fill] (1,1) circle (0.1);
\draw[fill] (2,1) circle (0.1);
\draw[fill] (3,1) circle (0.1);
\draw[fill] (4,1) circle (0.1);

\end{tikzpicture}
\end{center}

\caption{A geometric representation of $U_{2,4}$.}
\label{u24}
\end{figure}
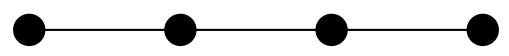

\begin{figure}[ht]
\begin{center}
\begin{tikzpicture}

%MK4 
\draw (11,1) -- (10,2.5) -- (9,4) -- (8,2.5) -- (7,1) ;
\draw (7,1) -- (10,2.5);
\draw (8,2.5) -- (11,1);
%\draw (9,1) -- (9,4);

\draw[fill] (7,1) circle (0.1);
%\draw[fill] (9,1) circle (0.1);
\draw[fill] (11,1) circle (0.1);
\draw[fill] (10,2.5) circle (0.1);
\draw[fill] (8,2.5) circle (0.1);
\draw[fill] (9,4) circle (0.1);
\draw[fill] (9,2) circle (0.1);

\end{tikzpicture}
\end{center}

\caption{The cycle matroid $M(K_4)$.}
\label{k4  matroid}
\end{figure}
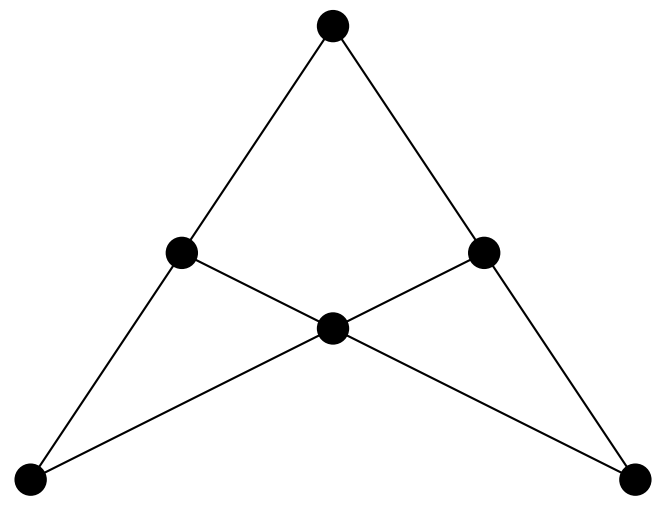

The main result of this paper, stated below, is proved in Section \ref{S.main}. Recall that $M(G)$ is the cycle matroid of a graph $G$.

\begin{restatable}{theorem}{outerplanar}
\label{T.outerplanar}
Let $G$ be a simple outerplanar graph. Then $M(G)$ is SZF-realizable. 
\end{restatable}

\section{Preliminaries}
\label{S.prelims}

All SZF-realizable uniform matroids are currently known, as stated in the following theorem (private communication).

\begin{theorem}
    \label{T.uniform}
    Let $r$ and $n$ be integers such that $0\leq r \leq n$. The matroid $U_{r,n}$ is SZF-realizable if and only if $r \in \{0, 1, n-2, n-1, n\}$ or $n\leq 6$.
\end{theorem}

\begin{corollary}
    \label{C.duality}
    The class of SZF-realizable matroids is not closed under duality.
\end{corollary}

In this section, we establish several fundamental facts regarding how the neighborhoods in an SZF-matroidal graph must relate to its corresponding matroid. The following well-known theorem (see, e.g. \cite{welsh} page 21 Theorem 1) will be used in the corollary that follows.

\begin{theorem}
\label{T.closure-ccts} 
Let $M$ be a matroid with $A\subseteq E(M)$. An element $x\in \cl(A)$ if and only if $x\in A$ or there exists a circuit $C$ of $M$ for which $C\setminus A=\{x\}$.
\end{theorem}

\begingroup
Let $G=(V,E)$ be a graph. We do not assume that $G$ is simple, but we will show that it is sufficient to consider graphs with at most one edge between pairs of distinct vertices and at most one loop at each vertex. For $u,v\in V$, let $\mu_G(u,v)$ be the edge multiplicity between $u$ and $v$. Then $\mu_G$ is symmetric and, if a vertex $u$ has no loops, then $\mu_G(u,u)=0$. If $\mu_G(u,v)=1$, we call the edge connecting $u$ and $v$ a \textit{single edge}. Set $N_G(u)=\{v\in V:\mu_G(u,v)>0\}$. Let $A\subseteq V$. In terms of skew zero forcing, if $A$ is the set of all currently filled vertices, then an element $u\in V(G)$ will force $v\in V(G)\setminus A$ exactly when \[\sum_{w\in V(G)\setminus A}\mu_G(u,w)=\mu_G(u,v)=1.\]

\begin{lemma}\label{outsidesequencelemma}
    Let $G$ be a finite multigraph, let $S\subseteq V(G)$ and let $t\in\textrm{scl}_G(S)\setminus S$. Then there exists a set $P\subseteq V(G)$ such that $t$ has a unique neighbor in $P$, joined to $t$ by a single edge. With $V(G)\setminus P$ initially filled, no vertex outside $S\cup\{t\}$ can force any other vertex. 
\end{lemma}

\begin{proof}
    For any initially filled set $S$ and vertex $x\in\textrm{scl}_G(S)$, we can construct a forcing sequence $S=S_0\subsetneq S_1\subsetneq\cdots\subsetneq S_k$, where $x\in S_k$ and, for each $i=1,\ldots,k$,
    \[S_i=S_{i-1}\cup\{x_i\},\quad N_G(v_i)\setminus S_{i-1}=\{x_i\},\quad\text{and}\quad \mu_G(v_i,x_i)=1\]
    for some vertices $x_i$ and $v_i$. Suppose we have such a sequence with minimum $k$ so that $t\in S_k$. Then $x_k=t$ and each $x_i$ is distinct and lies outside $S$. Each $v_i$ is also distinct since after it forces, all its neighbors are filled and it cannot force again. Set $P=\{v_1,\ldots,v_k\}$. Since $k$ is minimum, it follows that, for $i<k$, the vertex $t$ is still unfilled when $v_i$ forces $x_i$. Thus $v_i$ is not adjacent to $t$ and $v_k$ has a single adjacency to $t$. Setting $p=v_k$, it follows that $N_G(t)\cap P=\{p\}$ and $\mu_G(t,p)=1$, as desired.

    For every $i<k$, it follows that $x_i$ is adjacent to at least two members of $P$: the vertex $v_i$ that forced it and some $v_j$ with $j>i$. If not, then every $v_j$ with $j>i$ can still perform its forcing without $x_i$ filled and thus the sequence element $S_i$ can be omitted, contradicting that $k$ is minimum. Thus, if exactly $P$ is unfilled, then $x_i$ cannot force any other vertex.

    Finally, if $z\not\in S\cup \{x_1,\ldots,x_k=t\}$, then $z$ remains unfilled throughout the sequence. It cannot be adjacent to any $v_i$ since it would be an additional unfilled neighbor when $v_i$ forces $x_i$. Hence $N_G(z)\cap P=\varnothing$, and $z$ cannot force when exactly $P$ is unfilled. The conclusion follows.
\end{proof}

\begin{lemma}\label{parallelsym}
    Let $G$ be an SZF-matroidal multigraph with $u,v\in V(G)$. If $u$ and $v$ are adjacent, then $v\in\textrm{scl}(N_G(u)\setminus\{v\})$ if and only if $u\in\textrm{scl}(N_G(v)\setminus \{u\})$.
\end{lemma}
\begin{proof}
    If $u=v$, then the conclusion is clear and we thus assume $u\neq v$. Suppose, toward contradiction, that $v\not\in\textrm{scl}(N_G(u)\setminus\{v\})$ but $u\in\textrm{scl}(N_G(v)\setminus \{u\})$. Let $x$ be any neighbor of $u$ such that $\mu_G(u,x)=1$. Initially filling $N_G(u)\setminus\{x\}$ results in $u$ having exactly one unfilled neighbor, namely $x$, and therefore $x\in\textrm{scl}(N_G(u)\setminus \{x\})$. If $x=v$, then we have a contradiction and we thus assume $\mu_G(u,v)\geq 2$. Consequently, $x\neq v$ and $x\in\textrm{scl}(N_G(u)\setminus\{x\})$.

   Since $G$ is SZF-matroidal, Theorem \ref{T.closure-ccts} applied to $SZF(G)$ shows that there is a matroid circuit $C\subseteq N_G(u)$ containing $x$. Since $v\not\in \textrm{scl}(N_G(u)\setminus \{v\})$ it follows that $v\not\in C$ and thus $x\in\textrm{scl}(N_G(u)\setminus\{v,x\})$.

    Now, starting with only $N_G(u)\setminus\{v,x\}$ filled, the vertex $v$ can never become filled since $v\not\in\textrm{scl}(N_G(u)\setminus\{v\})$ implies $v\not\in\textrm{scl}(N_G(u)\setminus\{v,x\})$. While $v$ is unfilled, no neighbor of $v$ can force any other vertex and no vertex can force $v$ because $v$ lies outside the closure. Thus a forcing sequence using forcing vertices entirely outside of $N_G(v)$ is needed to establish that $x\in\textrm{scl}(N_G(u)\setminus\{v,x\})$.

    Set $S=N_G(v)\setminus\{u\}$ and $t=u$. Then Lemma \ref{outsidesequencelemma} applies and we can let $P\subseteq V(G)$ and $p\in P$ such that $N_G(u)\cap P=\{p\}${, and} $\mu_G(u,p)=1$, and, with $V(G)\setminus P$ initially filled, no vertex outside of $S\cup\{t\}=N_G(v)$ can force any other vertex. Setting $x=p$ in the previous paragraph gives an immediate contradiction. Using symmetry of $u$ and $v$, the conclusion follows.
\end{proof}

Following the notation of \cite{FreijHollantiEtAl2021}, we write $\textrm{cyc}_M$ for the cyclic operator of $M$:
\[\textrm{cyc}_M(S) = \{x \in S : r_M(S\setminus\{x\}) = r_M(S)\}.\] 
From closure definitions and Theorem \ref{T.closure-ccts}, we have the following.

\begin{corollary}\label{cyccor}
    If $M$ is a matroid and $S\subseteq E(M)$, then
    \[
    \textrm{cyc}_M(S) = \{x\in S : x\in\textrm{cl}_M(S\setminus\{x\})\} = \bigcup\{C : C \textrm{ is a circuit of } M \textrm{ and } C\subseteq S\}.
    \]
\end{corollary}

\begin{proposition}
\label{P.sortofsimple}
    If $M$ is SZF-realizable, then $M$ is SZF-realizable by a graph having at most one edge between each pair of distinct vertices and at most one loop at each vertex.
\end{proposition}
\begin{proof}
    Let $M$ have ground set $V$ and suppose $G$ is a graph with vertex set $V$ satisfying $\textrm{scl}_G(X)=\textrm{cl}_M(X)$ for all $X\subseteq V$. For each $u\in V$, set $C_u=\textrm{cyc}_M(N_G(u))$. Create a graph $H$ with vertices $V(G)$ and edge set $E_H=\{uv : v\in C_u\}$. Applying Lemma \ref{parallelsym} and Corollary \ref{cyccor}, we have that $v\in C_u$ if and only if $u\in C_v$. Thus $H$ is an undirected graph with no parallel edges and satisfies $N_H(u)=C_u$ for all $u\in V(H)$. Since $C_u\subseteq N_G(u)$, no new adjacency is introduced. If $\mu_G(u,v)=1$, then $u$ can force $v$ from $N_G(u)\setminus\{v\}$, so $v\in C_u$. Thus every original single edge is retained.

    Choose $A\subseteq V$. We first show $\textrm{scl}_G(A)\subseteq\textrm{scl}_H(A)$. Suppose vertex $u$ forces vertex $v$ in $G$. Since forcing from $u$ to $v$ is only possible when $\mu_G(u,v)=1$, every possible force from $u$ is retained in $H$. For the reverse inclusion, consider a force from $u$ to $v$ in $H$. Then $v\in C_u$ and there is a matroid circuit $C\subseteq N_G(u)$ containing $v$. Every element of $C$ belongs to $C_u$, and thus $C\subseteq C_u$. Therefore $v\in\textrm{cl}_M(C\setminus\{v\})\subseteq\textrm{cl}_M(C_u\setminus\{v\})$. In order for $u$ to force $v$ in $H$, all of $C_u\setminus\{v\}$ must be filled. Continuing in this way for the entire forcing sequence, we have that every element of $\textrm{scl}_H(A)$ is contained in $\textrm{cl}_M(A)$. By the construction of $G$, then, $\textrm{scl}_H(A)\subseteq \textrm{cl}_M(A) = \textrm{scl}_G(A)$ and we can conclude that $\textrm{scl}_G(A)=\textrm{scl}_H(A)$. Since $H$ has only single edges, the conclusion follows.
\end{proof}

Unless otherwise stated, all graphs considered from this point onward have at most one edge between each pair of distinct vertices and at most one loop at each vertex. The above result shows that this does not reduce the size of the class of SZF-realizable matroids, though this convention significantly reduces complexity when studying it.

\endgroup

\begin{proposition}
\label{P.seriesneighborhoods}
    Let $G$ be an SZF-matroidal graph with matroid $M$ and take $x,y \in E(M)$, where $x\neq y$, and neither $x$ nor $y$ are coloops of $M$. Then $n(x)=n(y)$ in $G$ if and only if $x$ and $y$ are in series in $M$.
\end{proposition}
\begin{proof}
    Assume $n(x)=n(y)\neq \emptyset$ in $G$, and consider the corank of $\{x,y\}$ in $M$.
    $$\begin{array}{r l}
         r^*_M(\{x,y\}) = & r(E-\{x,y\}) + |\{x,y\}| - r(M) \\
         = & r(E-\{x,y\}) + 2 - r(M) 
    \end{array} $$
    Observe that the set $E-\{x,y\}$ is closed with respect to SZF in $G$. To see this, fill $E-\{x,y\}$ and note that every neighbor of $x$ and $y$ has two unfilled neighbors, namely $x$ and $y$. The neighborhoods of $x$ and $y$ are equal and either they are both completely filled (if $\{x,y\} \not \subseteq n(x)=n(y)$) or they have two unfilled vertices (if $\{x,y\} \subseteq n(x)=n(y)$), and thus will not force any new vertices. Every other vertex in $G$ has a neighborhood that is completely filled and will not force new vertices. Thus, $\cl(E-\{x,y\})=E-\{x,y\}$. Now consider the set $(E-\{x,y\})\cup y = E-\{x\}$ and note that $\cl(E-\{x\})= E$. This implies that $r(E-\{x\}) = r(M)$ and $r(E-\{x,y\}) = r(M)-1$. So, 
    $$\begin{array}{r l}
         r^*_M(\{x,y\}) = & r(E-\{x,y\}) + 2 - r(M) \\
         = & r(M)-1+2-r(M) \\
         = & 1 .\\
    \end{array} $$
    Thus, as neither $x$ nor $y$ is a coloop, $x$ and $y$ must be in series.

    Now assume that $x$ and $y$ are in series in $M$. Then $r^*_M(\{x,y\})=1$ and thus $r(E-\{x,y\}) = r(M) -1$. Moreover, neither $x$ nor $y$ is contained in $\cl(E-\{x,y\})$, since any circuit of $M$ that contains $x$ must also contain $y$ and vice versa. Assume toward contradiction that there is some $a\in n(x)-n(y)$. Consider $n(a)$ and $\cl(E-\{x,y\})$. Every vertex in $n(a)$ is filled by $\cl(E-\{x,y\})$ except $x$. This implies that $x$ should be forced in $\cl(E-\{x,y\})$, a contradiction. A symmetric argument shows there can be no elements in $n(y)-n(x)$, and thus $n(x)=n(y)$. 
    
\end{proof}

\begin{definition} {Let $G$ be a graph and let $M$ be a matroid with ground set $V(G)$. The} neighborhood of {a vertex in} $G$ is \textit{cyclic in $M$} if it is a union of circuits of $M$.
\end{definition}

\begin{proposition}
\label{P.closurematch}
{Let $G$ be a graph and let $M$ be a matroid with ground set $V(G)$.} Suppose that every circuit of $M$ is a neighborhood of $G$ and every neighborhood of $G$ is cyclic in $M$. Then, for every $S\subseteq V(G)$, we get $\textrm{scl}_G(S)=\cl_M(S)$.
\end{proposition}
\begin{proof}
    Let $S=S_0\subsetneq S_1\subsetneq\cdots\subsetneq S_k=\textrm{scl}_G(S)$ be a forcing sequence, where for each $i=1,\ldots,k$,
    \[S_i=S_{i-1}\cup\{x_i\}\quad\text{and}\quad N(v_i)\setminus S_{i-1}=\{x_i\}\]
    for some vertex $v_i$. We prove by induction on $i$ that $S_i\subseteq\cl_M(S)$. The case $i=0$ is immediate. Suppose that $S_{i-1}\subseteq\cl_M(S)$. Since $N(v_i)$ is cyclic and contains $x_i$, there is a circuit $C$ such that $x_i\in C\subseteq N(v_i)\subseteq S_{i-1}\cup\{x_i\}$. Then, by Theorem \ref{T.closure-ccts}, we have that $x_i\in\cl_M(S_{i-1})\subseteq\cl_M(S)$. Hence $S_i\subseteq\cl_M(S)$, completing the induction.

    Conversely, suppose that $x\in\cl_M(S)\setminus S$. By Theorem \ref{T.closure-ccts}, there is a circuit $C$ with $x\in C\subseteq S\cup\{x\}$. By our hypothesis, let $v\in V(G)$ so that $N(v)=C$. Then $N(v)\setminus S=\{x\}$ and $v$ forces $x$. Thus $\cl_M(S)\subseteq \textrm{scl}_G(S)$, and the result holds.

\end{proof}

\begin{proposition}
\label{P.circuitflats}
    Let $M$ be an SZF-realizable matroid. Suppose $C$ is a circuit-flat of $M$, and let $G$ be a graph such that $SZF(G)\cong M$. Then there is a vertex $v \in V(G)$ such that $n(v) = C$. 
\end{proposition}
\begin{proof}
    This follows immediately from the observation that since $C\backslash e$ is independent for all $e\in C$ and $C$ is a flat, we know that filling $C\ba e$ in $G$ will force $e$ and no other vertices of $G$.
\end{proof}

\begin{corollary}
\label{C.circuitflats}
    Let $M$ be a matroid. If the number of circuit-flats of $M$ exceeds the number of elements in $M$, then $M$ is not SZF-realizable.
\end{corollary}

The previous corollary and the proposition that follows will both be used extensively in Section \ref{S.main}. 

\begin{proposition}
\label{P.symmetricdifference}
    Let $\mathcal{N}$ be a set of neighborhoods for a graph $G$, and assume there are neighborhoods $N_1,N_2 \in \mathcal{N}$ such that $N_1\cap N_2 = \{e\}$. Adding the set $K=N_1 \triangle N_2$ to the set $\mathcal{N}$ does not change $\textrm{scl}_G(S)$ for any $S\subseteq V(G)$.
\end{proposition}
\begin{proof}
    Observe that the set $K$ will only force vertices in the situation where $K\backslash f$ {is filled in $G$}, for some $f\in K$. If $f=e$, then $N_1$ forces $e$ and $K$ has no effect on $\textrm{scl}_G$. So, without loss of generality, assume $f\in N_1-N_2$. Note that $N_2\ba e \subseteq K\ba f$, so filling $K\ba f$ will result in $e$ being forced because of $N_2$, which in turn will force $f$ because of $N_1$. Thus, adding the set $K$ to $\mathcal{N}$ does not change how the SZF closure operator works in the graph. 
\end{proof}

\section{Minors}
\label{S.minors}

As mentioned in Sections \ref{S.intro} and \ref{S.prelims}, our research considers skew zero forcing (SZF) on multigraphs. The primary reason for this is that allowing loops and multiple edges makes graphs closed under the minor operation. There are many seminal results in combinatorics related to describing classes of minor-closed graphs and matroids via their excluded minors. Another good reason to consider graph loops is that they have an interesting relationship to their corresponding matroid elements when a graph is SZF-matroidal. 

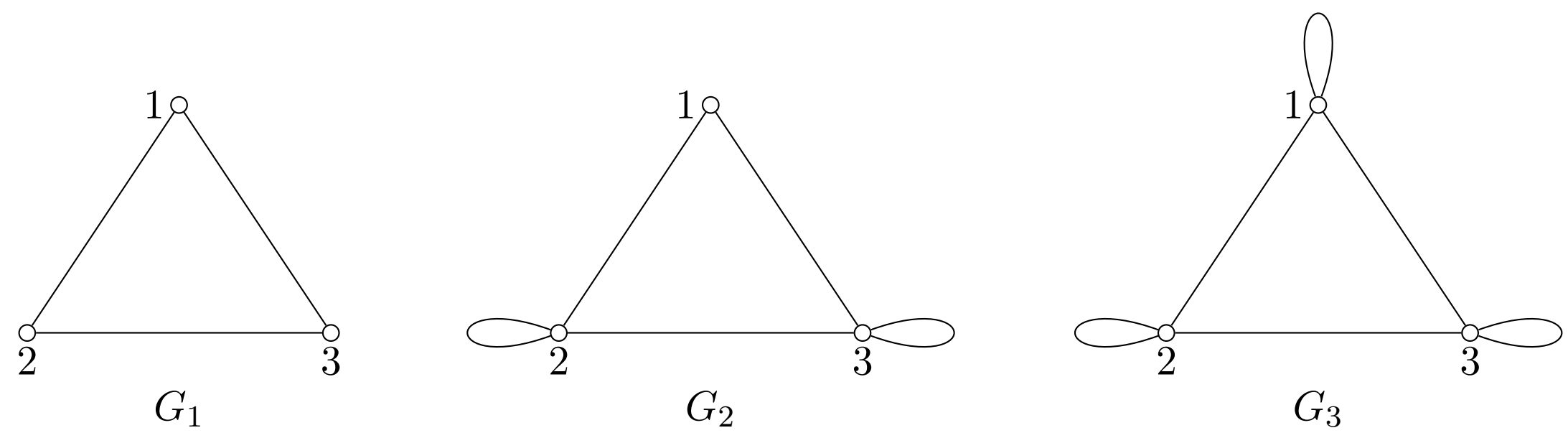
\begin{figure}[ht]

\begin{center}
\begin{tikzpicture}[scale=1.4,
  vertex/.style={circle,draw,fill=white,inner sep=1.5pt}
]

%
% no loops
%
\node[vertex] (v2) at  (0.5,7) {};
\node[vertex] (v3) at  (2.5,7) {};
\node[vertex] (v1) at  (1.5,8.5) {};

\node[left] at (v1) {1};
\node[below] at (v2) {2};
\node[below] at (v3) {3};

\draw (v1) -- (v2) -- (v3) -- (v1);

\node at (1.5,6.5) {$G_1$};

%
% 2 loops
%
\node[vertex] (v4) at  (4,7) {};
\node[vertex] (v5) at  (5,8.5) {};
\node[vertex] (v6) at  (6,7) {};

\node[below] at (v4) {2};
\node[left] at (v5) {1};
\node[below] at (v6) {3};

\draw (v4) -- (v5) -- (v6) -- (v4);

\draw (v4) to[out=160, in=200, looseness=50] (v4);
\draw (v6) to[out=20, in=340, looseness=50] (v6);

\node at (5,6.5) {$G_2$};

%
% 3 loops
%
\node[vertex] (v7) at  (8,7) {};
\node[vertex] (v8) at  (9,8.5) {};
\node[vertex] (v9) at  (10,7) {};

\node[below] at (v7) {2};
\node[left] at (v8) {1};
\node[below] at (v9) {3};

\draw (v7) -- (v8) -- (v9) -- (v7);

\draw (v7) to[out=160, in=200, looseness=50] (v7);
\draw (v9) to[out=20, in=340, looseness=50] (v9);
\draw (v8) to[out=70, in=110, looseness=50] (v8);

\node at (9,6.5) {$G_3$};

\end{tikzpicture}
\end{center}

\caption{Three graphs $G_1$, $G_2$, and $G_3$.}
\label{F.loops}

\end{figure}

\begin{figure}[ht]

\begin{center}
\begin{tikzpicture}[scale=1.4,
  vertex/.style={circle,draw,fill=white,inner sep=1.5pt}
]

%
% no loops
%
\node[vertex] (v2) at  (0.5,7) {};
\node[vertex] (v3) at  (2.5,7) {};
\node[vertex] (v1) at  (1.5,8.5) {};

\node[left] at (v1) {1};
\node[below] at (v2) {2};
\node[below] at (v3) {3};

\draw (v1) -- (v2) -- (v3) -- (v1);

\node at (1.5,6.5) {$H_1$};

%
% 2 loops
%
\node[vertex] (v4) at  (4,7) {};
\node[vertex] (v5) at  (5,8.5) {};
\node[vertex] (v6) at  (6,7) {};

\node[below] at (v4) {2};
\node[left] at (v5) {1};
\node[below] at (v6) {3};

\draw (v4) -- (v5) -- (v6) -- (v4);

\draw (v4) to[out=30, in=150, looseness=1] (v6);

\node at (5,6.5) {$H_2$};

\end{tikzpicture}
\end{center}

\caption{Two graphs $H_1$ and $H_2$.}
\label{F.parallel}

\end{figure}
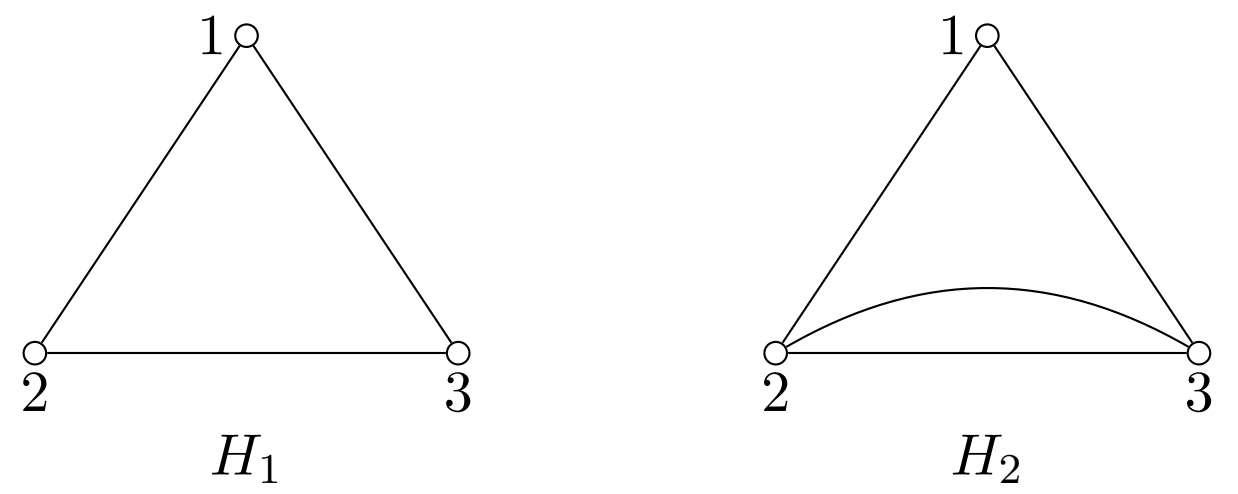

The addition or removal of a loop or parallel edge can change a graph from SZF-matroidal to non-SZF-matroidal and vice versa. For example, consider the graphs $G_1$, $G_2$, and $G_3$ depicted in Figure \ref{F.loops}. The graph $G_1$ is SZF-matroidal with $SZF(G_1)\cong U_{1,3}$. The graph $G_3$ is SZF-matroidal with $SZF(G_3)\cong U_{2,3}$. However, the graph $G_2$ is not SZF-matroidal. To see this, observe that $\textrm{scl}_{G_2}(\{2\})=\{1,2,3\}$, but $\textrm{scl}_{G_2}(\{1\})=\{1\}$, which violates the closure exchange axiom for matroids. A similar situation exists for parallel edges, where the graph $H_1$ in Figure \ref{F.parallel} is SZF-matroidal with $SZF(H_1)\cong U_{1,3}$ but $H_2$ is not SZF-matroidal.

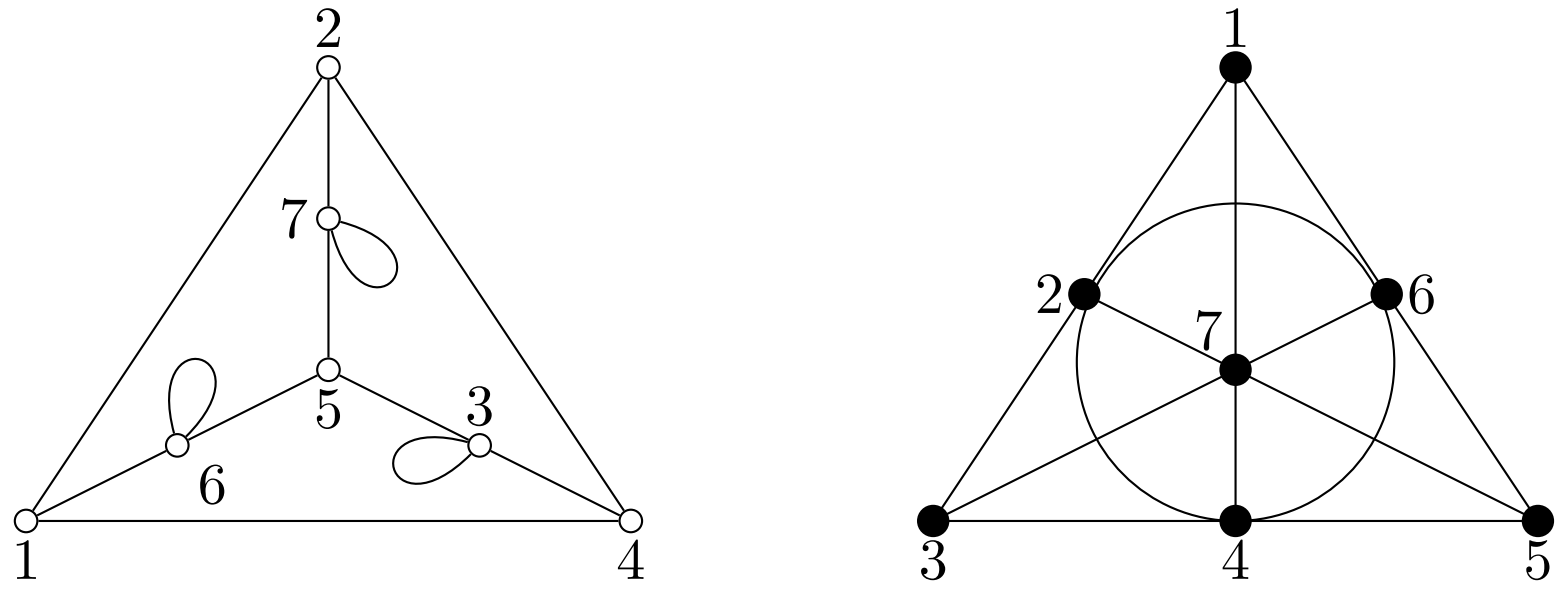
\begin{figure}[ht]
\begin{center}
\begin{tikzpicture}[scale=1,
  vertex/.style={circle,draw,fill=white,inner sep=1.5pt}
]

\begin{scope}[shift={(-6,0)}]
    \node[vertex] (a) at (7,1) {};
    \node[vertex] (b) at (8, 1.5) {};
    \node[vertex] (c) at (11,1) {};
    \node[vertex] (d) at (10,1.5) {};
    \node[vertex] (e) at (9,3) {};
    \node[vertex] (f) at (9,4) {};
    \node[vertex] (g) at (9, 2) {};

    \node[above] at (f) {2};
    \node[below right] at (b) {6};
    \node[below] at (c) {4};
    \node[left] at (e) {7};
    \node[below] at (a) {1};
    \node[above] at (d) {3};
    \node[below] at (g) {5};

    \draw (a) -- (b) -- (g) -- (e) -- (f) -- (a) -- (c) -- (d) -- (g);
    \draw (c) -- (f);
    
    \myloop{b}{45};
    \myloop{d}{165};
    \myloop{e}{-75};

\end{scope}

\if{false}
%fano graph
\begin{scope}
\node[vertex] (v6) at  (1,1) {};
\node[vertex] (v5) at  (3,1) {};
\node[vertex] (v3) at  (5,1) {};
\node[vertex] (v4) at  (4,4) {};
\node[vertex] (v1) at  (2,4) {};
\node[vertex] (v2) at  (3,3) {};
\node[vertex] (v7) at  (3,2) {};

\node [below right] at  (v6) {6};
\node [below] at  (v5) {5};
\node [below left] at  (v3) {3};
\node [right] at  (v4) {4};
\node [left] at  (v1) {1};
\node [right] at  (v2) {2};
\node [right] at  (v7) {7};

\draw (v6) -- (v5) -- (v3) -- (v4) -- (v1) -- (v6) ;
\draw (v5) -- (v7) -- (v2) -- (v1);
\draw (v2) -- (v4);

\draw (v7) to[out=160, in=200, looseness=50] (v7);
\draw (v6) to[out=160, in=200, looseness=50] (v6);
\draw (v3) to[out=20, in=340, looseness=50] (v3);
\end{scope}
\fi

%fano matroid
\begin{scope}
\draw (7,1) -- (9,1) -- (11,1) -- (10,2.5) -- (9,4) -- (8,2.5) -- (7,1) ;
\draw (7,1) -- (10,2.5);
\draw (8,2.5) -- (11,1);
\draw (9,1) -- (9,4);
\draw(9,2.05) circle (1.05);

\draw[fill] (7,1) circle (0.1);
\draw[fill] (9,1) circle (0.1);
\draw[fill] (11,1) circle (0.1);
\draw[fill] (10,2.5) circle (0.1);
\draw[fill] (8,2.5) circle (0.1);
\draw[fill] (9,4) circle (0.1);
\draw[fill] (9,2) circle (0.1);

\node [anchor=north] at  (7,1) {3};
\node [anchor=north] at  (9,1) {4};
\node [anchor=north] at  (11,1) {5};
\node [anchor=west] at  (10,2.5) {6};
\node [anchor=east] at  (8,2.5) {2};
\node [anchor=south] at  (9,4) {1};
\node [anchor=305] at  (9,2) {7};
\end{scope}

\end{tikzpicture}
\end{center}

\caption{A graph whose SZF closure gives rise to the Fano matroid, $F_7$, and a geometric representation of $F_7$.}
\label{F.fano}
\end{figure}

Moreover, some matroids that are SZF-realizable require graph loops or parallel edges in their SZF-matroidal graphs. One example of this is the Fano matroid, $F_7$, depicted on the right in Figure \ref{F.fano}. A graph whose SZF closure gives rise to the Fano matroid is depicted on the left in Figure \ref{F.fano}, meaning that the Fano matroid is SZF-realizable. However, every graph that gives rise to the Fano matroid must have loops or parallel edges, as noted in the following proposition. 

\begin{proposition}
Let $G$ be a graph such that $SZF(G)\cong F_7$. Then $G$ must contain a loop or parallel edge.
\begin{proof}
We have already established that $F_7$ is SZF-realizable, so let $G$ be a graph in $F_{7_{SZF}}$. Observe that $F_7$ has seven distinct circuit-flats. By Proposition \ref{P.circuitflats}, we know that each circuit-flat must correspond exactly to a neighborhood of $G$. Since $F_7$ has seven elements, we know that $G$ must have seven vertices and thus each neighborhood of $G$ must have size exactly three. Assume $G$ has no loops or parallel edges. Then the size of each neighborhood corresponds exactly to the degree of its vertex. The First Theorem of Graph Theory states that the sum of all degrees in a graph equals twice the number of edges. However, the sum of all degrees of $G$ would be 21 in this case, an odd number. This contradiction implies that $G$ must have a loop or parallel edge.
\end{proof}
\label{P.fano has loops}
\end{proposition}

%
% Fano graph and the matroid side by side
%

We now consider the minor operation for both graphs and matroids with respect to SZF, beginning with matroids. There are circumstances where deletion and contraction maintain the SZF-realizability of a matroid (matroid coloops, for example), but the class of SZF-realizable matroids is not minor-closed. In fact, the following two theorems establish that neither SZF-realizable matroids nor SZF-matroidal graphs form a minor-closed class.

\begin{proposition}\label{sixelement-notmatroidal}
    The simple rank-three matroid on six elements consisting of a four-point line and two elements placed freely is not SZF-realizable.
\end{proposition}
\begin{proof}
    Let $M$ be the matroid described in the hypothesis where $L=\{a,b,c,d\}$ is the four-point line and $F=\{y,z\}$ are the remaining points. Suppose $G$ realizes $M$. Since $L$ is a hyperplane, $\{y,z\}$ is a cocircuit. As neither $y$ nor $z$ is a coloop, Proposition \ref{P.seriesneighborhoods} requires that $N_G(y)=N_G(z)$. Since $M$ is simple and has no coloops, every neighborhood of $G$ has at least three vertices. Every neighborhood of size exactly three is dependent and is thus contained in $L$. Further, every pair of $L$ must be contained in at least one size-three neighborhood. At least three distinct such neighborhoods are necessary in $G$, since any fewer would be not able to cover each of the six pairs in $L$.

    Now, if $N_G(v)$ is one of these {size-three} neighborhoods, then $y\not\in N_G(v)$ and symmetry requires that $v\not\in N_G(y)$. Then at least three vertices lie outside $N_G(y)$ and thus $|N(y)|\leq 3$ which implies that $|N_G(y)|=3$. Then $N_G(y)\subseteq L$. However, any element in $N_G(y)$ must have a neighborhood containing $y$. Thus only those elements in $V(G)\setminus N_G(y)$, namely $y$, $z$, and one member of $L$ can have neighborhoods contained in $L$. Since $N_G(y)=N_G(z)$, this contradicts the need for three distinct such neighborhoods. The conclusion follows.
\end{proof}

\begin{theorem}
The class of SZF-realizable matroids is not closed with respect to the minor operation.
\label{T.minor closed matroid}
\end{theorem}

\begin{figure}[ht]
\begin{center}
\begin{tikzpicture}[scale=1,
  vertex/.style={circle,draw,fill=white,inner sep=1.5pt}
]

% Graph coordinates
\node[vertex] (T)  at (3,4) {};
\node[vertex] (UL) at (1,2.5) {};
\node[vertex] (C)  at (3,2) {};
\node[vertex] (UR) at (5,2.5) {};
\node[vertex] (BL) at (1,1) {};
\node[vertex] (B)  at (3,0.5) {};
\node[vertex] (BR) at (5,1) {};

% Graph edges
\draw
    (T) -- (UL)
    (T) -- (C)
    (T) -- (UR)
    (UL) -- (C)
    (C) -- (UR)
    (UL) -- (BL)
    (C) -- (B)
    (UR) -- (BR)
    (BL) -- (B)
    (B) -- (BR);

% Loops
\draw (T) to[loop above, min distance=1.5cm] (T);
\draw (BL) to[loop left, min distance=1.5cm] (BL);
\draw (BR) to[loop right, min distance=1.5cm] (BR);

\coordinate[label=above right:{$b$}] (T) at (T);
\coordinate[label=above right:{$c$}] (T) at (C);
\coordinate[label=above right:{$z$}] (T) at (UR);
\coordinate[label=above left:{$a$}] (T) at (BL);
\coordinate[label=above left:{$x$}] (T) at (B);
\coordinate[label=above right:{$d$}] (T) at (BR);
\coordinate[label=above left:{$y$}] (T) at (UL);

% Matroid 
\draw (11,1) -- (10,2.5) -- (9,4) -- (8,2.5) -- (7,1) ;
\draw (7,1) -- (11,1);

\draw[fill] (7,1) circle (0.1);
\draw[fill] (11,1) circle (0.1);
\draw[fill] (10,2.5) circle (0.1);
\draw[fill] (8,2.5) circle (0.1);
\draw[fill] (9,4) circle (0.1);
\draw[fill] (8.25,1) circle (0.1);
\draw[fill] (9.75,1) circle (0.1);

\coordinate[label=above right:{$d$}] (T) at (11,1);
\coordinate[label=above right:{$z$}] (T) at (10,2.5);
\coordinate[label=above right:{$x$}] (T) at (9,4);
\coordinate[label=above left:{$y$}] (T) at (8,2.5);
\coordinate[label=above left:{$a$}] (T) at (7,1);
\coordinate[label=above right:{$b$}] (T) at (8.25,1);
\coordinate[label=above right:{$c$}] (T) at (9.75,1);

\end{tikzpicture}
\end{center}

\caption{A graph and its corresponding SZF matroid.}
\label{F.notminorclosed}
\end{figure}
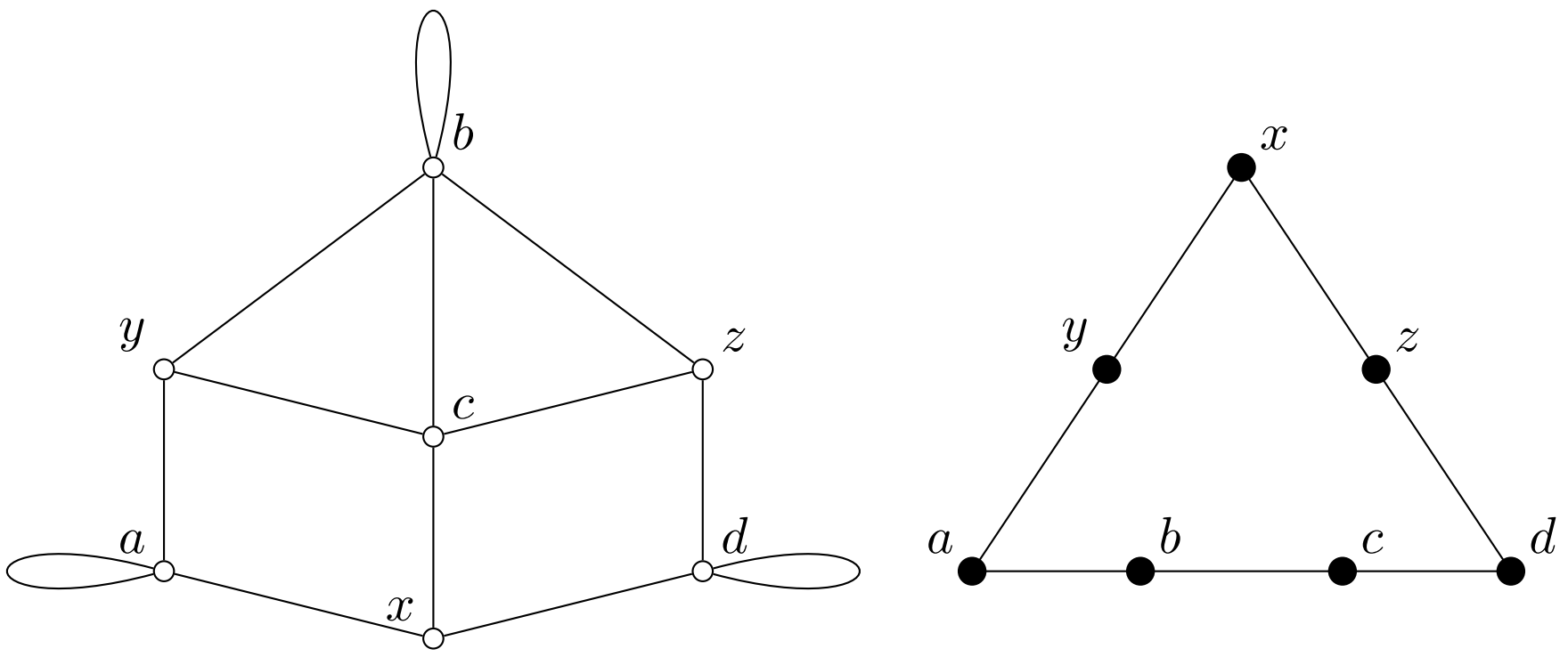

\begin{proof}
Consider the matroid on the right in Figure \ref{F.notminorclosed}. Using the graph on the left of Figure \ref{F.notminorclosed}, we observe that the matroid $M$ is SZF-realizable. However by Proposition \ref{sixelement-notmatroidal}, the matroid $M\backslash x$ is not.
\end{proof}

\begin{theorem}
The class of SZF-matroidal graphs is not closed with respect to vertex deletion, edge deletion, or edge contraction.
\label{T.minor closed graph}
\end{theorem}

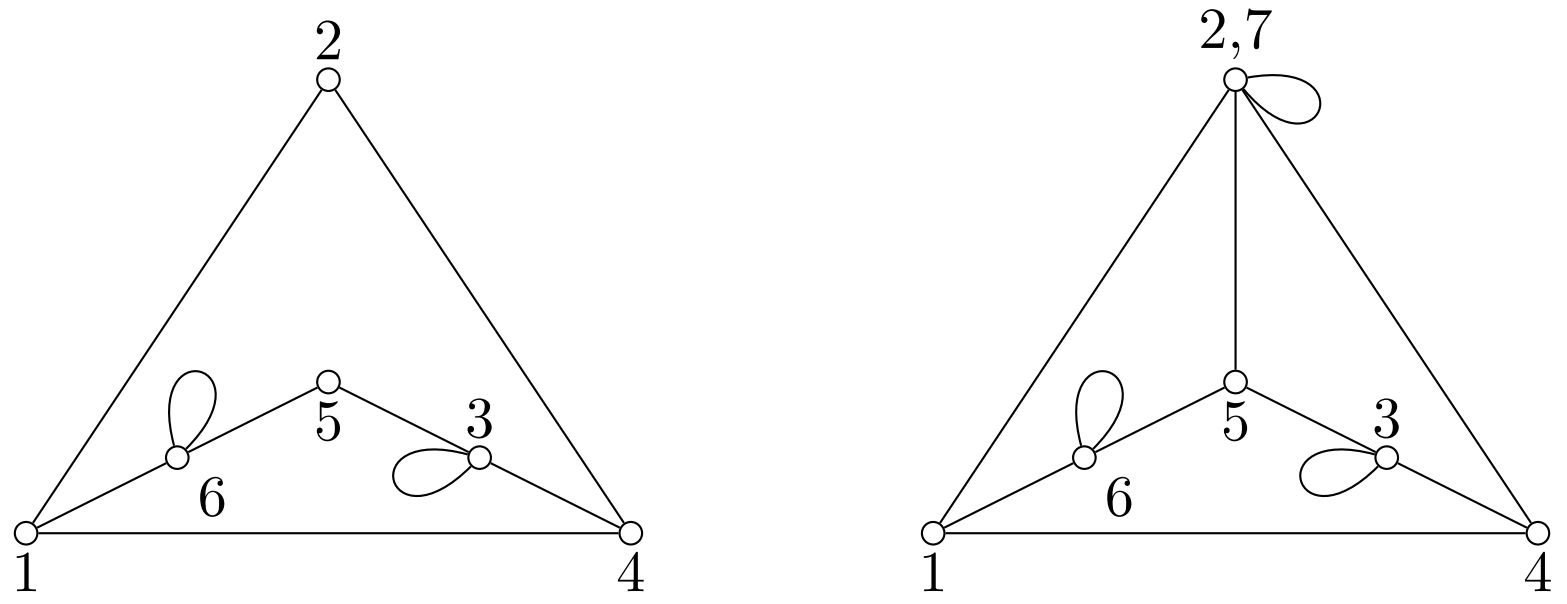
\begin{figure}[ht]
\begin{center}
\begin{tikzpicture}[scale=1,
  vertex/.style={circle,draw,fill=white,inner sep=1.5pt}
]

%DELETION in the Fano graph
\begin{scope}[shift={(-6,6)}]
    \node[vertex] (a) at (7,1) {};
    \node[vertex] (b) at (8, 1.5) {};
    \node[vertex] (c) at (11,1) {};
    \node[vertex] (d) at (10,1.5) {};
%   \node[vertex] (e) at (9,3) {};
    \node[vertex] (f) at (9,4) {};
    \node[vertex] (g) at (9, 2) {};

    \node[above] at (f) {2};
    \node[below right] at (b) {6};
    \node[below] at (c) {4};
%   \node[left] at (e) {7};
    \node[below] at (a) {1};
    \node[above] at (d) {3};
    \node[below] at (g) {5};

    \draw (a) -- (b) -- (g);
    \draw (f) -- (a) -- (c) -- (d) -- (g);
    \draw (c) -- (f);
    
    \myloop{b}{45};
    \myloop{d}{165};
\end{scope}

% CONTRACTION in the Fano graph
\begin{scope}[shift={(0,6)}]
    \node[vertex] (a) at (7,1) {};
    \node[vertex] (b) at (8, 1.5) {};
    \node[vertex] (c) at (11,1) {};
    \node[vertex] (d) at (10,1.5) {};
%    \node[vertex] (e) at (9,3) {};
    \node[vertex] (f) at (9,4) {};
    \node[vertex] (g) at (9, 2) {};

    \node[above] at (f) {2,7};
    \node[below right] at (b) {6};
    \node[below] at (c) {4};
%   \node[above right] at (f) {(7)};
    \node[below] at (a) {1};
    \node[above] at (d) {3};
    \node[below] at (g) {5};

    \draw (a) -- (b) -- (g);
    \draw (f) -- (a) -- (c) -- (d) -- (g);
    \draw (c) -- (f) -- (g);
    
    \myloop{b}{45};
    \myloop{d}{165};
    \myloop{f}{310};
\end{scope}

\end{tikzpicture}
\end{center}

\caption{Two graphs that do not give rise to a matroid via the SZF closure operator.}
\label{F.fano graph operations}
\end{figure}

\begin{proof}
We first address vertex deletion and consider the graph $G$ on the left in Figure \ref{F.fano}. We know that $SZF(G)\cong F_7$, with labels as in Figure \ref{F.fano}. Consider $G\backslash 7$, depicted on the left in Figure \ref{F.fano graph operations}. In $G\ba 7$ observe that $\textrm{scl}(\{2\})=\{2\}$, $\textrm{scl}(\{1,2\})=\{1,2,3,4,5,6\}$, but $\textrm{scl}(\{2,5\})=\{2,5\}$, which violates our exchange axiom. To see that the class is not closed under edge deletion, consider Figure \ref{F.loops} and observe that $G_2$ is obtained from $G_3$ by deleting an edge, but $G_2$ is not SZF-matroidal while $G_3$ is. Finally, for edge contraction, consider contracting the edge 27 in the graph $G$ depicted on the left in Figure \ref{F.fano}. This gives rise to the graph on the right in Figure \ref{F.fano graph operations}. As $\textrm{scl}(\{1,4\})=\{1,4\}$, but  $\textrm{scl}(\{1,2\})=\{1,2,3,4,5,6\}$, this graph is not SZF-matroidal.
\end{proof}

\begin{corollary}
The class of SZF-matroidal graphs is not closed with respect to the minor operation.
\label{C.minor closed graph}
\end{corollary}

\section{Graphs that give rise to skew zero forcing matroids}
\label{S.main}

In this section, we investigate the {SZF-realizability of graphic matroids}. In many of the arguments in this section, we begin with a graph $G$ and attempt to build an associated graph $M(G)_{SZF}$ that gives rise to $M(G)$ via the SZF closure operator. We will often attempt to construct $M(G)_{SZF}$ directly from $G$ by defining the neighborhoods of $M(G)_{SZF}$ using the notation $n(e)$ where $e$ is an edge of $G$, an element of $M(G)$, and a vertex of $M(G)_{SZF}$. As we will see in this section, there are graphic matroids that are not SZF-realizable. This section explores some common and important graph classes to determine which of their cycle matroids are SZF-realizable.

Observe that if $G$ is an acyclic graph on $n$ edges, then $M(G)\cong U_{n,n}$, which is SZF-realizable by Theorem \ref{T.uniform}. In fact, edges that do not appear in any cycles of a graph have no effect on the associated graphic matroid's SZF-realizability. 

\begin{proposition}
\label{P.2conn}
For a loopless graph $G$,  if the matroids associated with each $2$-connected block of $G$ are SZF-realizable, then $M(G)$ is SZF-realizable.
\end{proposition}

\begin{proof}
First note that if $G$ has any isolated vertices, they have no effect on $M(G)$ and can be ignored. If $G$ is disconnected, then $M(G)$ is disconnected and we may consider each component separately, constructing separate components of $M(G)_{SZF}$ for each connected component of $M(G)$. Thus, we may assume $G$ is connected. If $G$ has edges that are not contained in any cycle of $G$, these edges correspond to coloops in $M(G)$, which are realized in $M(G)_{SZF}$ as isolated vertices, since coloops are exactly the elements of $M(G)$ that do not show up in any circuit. Thus, we can assume that $G$ has no such edges. If $G$ is not $2$-connected but every edge of $G$ is in some cycle of $G$, then $G$ contains cut vertices that separate $2$-connected blocks of $G$. In this case, the matroid $M(G)$ is not connected and, similar to above, we may construct separate components of $M(G)_{SZF}$ for each $2$-connected block of $M(G)$.    
\end{proof}

We begin our investigation of graph classes by considering cycles, denoted $C_n$, complete graphs, denoted $K_n$, and complete bipartite graphs, denoted $K_{m,n}$.

\begin{proposition}
\label{P.cyclematroid}
The matroid $M(C_n)$ is SZF-realizable.
\end{proposition}

\begin{proof}
    This follows immediately from the fact that $M(C_n)\cong U_{n-1,n}$ and Theorem \ref{T.uniform}.
\end{proof}

\begin{proposition}
\label{P.completeno}
The matroid $M(K_n)$ is not SZF-realizable when $n\geq 6$.
\end{proposition}

\begin{proof}
Observe that $K_n$ has $\binom{n}{2}$ edges and $\binom{n}{3}$ triangles, each of which is a circuit-flat of $M(K_n)$. Thus, Corollary \ref{C.circuitflats} gives us the result when $n\geq 6$.
\end{proof}

\begin{proposition}
\label{P.completeyes}
The matroid $M(K_n)$ is SZF-realizable when $n\in\{1,2,3,4,5\}$.
\end{proposition}

\begin{figure}
\begin{tikzpicture}[scale=1.4,
  vertex/.style={circle,draw,fill=white,inner sep=1.5pt}
]

% curved edge
\draw plot [smooth] coordinates {(0,2)(2,3)(4,2)};

% vertices (top, middle, bottom, left, right)
\node[vertex] (ll) at (0,0) {};      % lower left
\node[vertex] (lm) at (2,0) {};      % lower middle
\node[vertex] (lr) at (4,0) {};     % lower right
\node[vertex] (ul) at (0,2) {};     % upper left
\node[vertex] (um) at (2,2) {};     % upper middle
\node[vertex] (ur) at (4,2) {};      % upper right

% straight edges
\draw (ll) -- (lm) -- (lr) -- (ur) -- (um) -- (ul) -- (ll);
\draw (lm) -- (um);

% loop edges
\draw (ll) to[out=290, in=250, looseness=50] (ll);
\draw (lm) to[out=290, in=250, looseness=50] (lm);
\draw (lr) to[out=290, in=250, looseness=50] (lr);
\draw (um) to[out=110, in=70, looseness=50] (um);

\end{tikzpicture} \hspace{0.5in}
\begin{tikzpicture}[scale=1.4,
  vertex/.style={circle,draw,fill=white,inner sep=1.5pt}
]

\def\Rout{2}
\def\Rin{1}
\def\rot{-126}  % rotation offset in degrees

% Outer 5-cycle (radius Rout), rotated so a flat edge is at the bottom
\foreach \k in {0,...,4} {
  \node[vertex] (o\k) at ({\rot + 72*\k}:\Rout) {};
}

% Inner 5-cycle (radius Rin), same rotation
\foreach \k in {0,...,4} {
  \node[vertex] (i\k) at ({\rot + 72*\k}:\Rin) {};
}

% Outer cycle edges
\foreach \k in {0,...,4} {
  \pgfmathtruncatemacro{\next}{mod(\k+1,5)}
  \draw (o\k) -- (o\next);
}

% Spokes (outer to inner)
\foreach \k in {0,...,4} {
  \draw (o\k) -- (i\k);
}

% Inner star edges: connect i_k to i_{k+2} (mod 5)
\foreach \k in {0,...,4} {
  \pgfmathtruncatemacro{\jump}{mod(\k+2,5)}
  \draw (i\k) -- (i\jump);
}

\end{tikzpicture}
    \caption{Graphs $M(K_4)_{SZF}$ and $M(K_5)_{SZF}$ that give rise to $M(K_4)$ and $M(K_5)$, respectively, via the SZF closure operator.}
    \label{F.k4k5}
\end{figure}
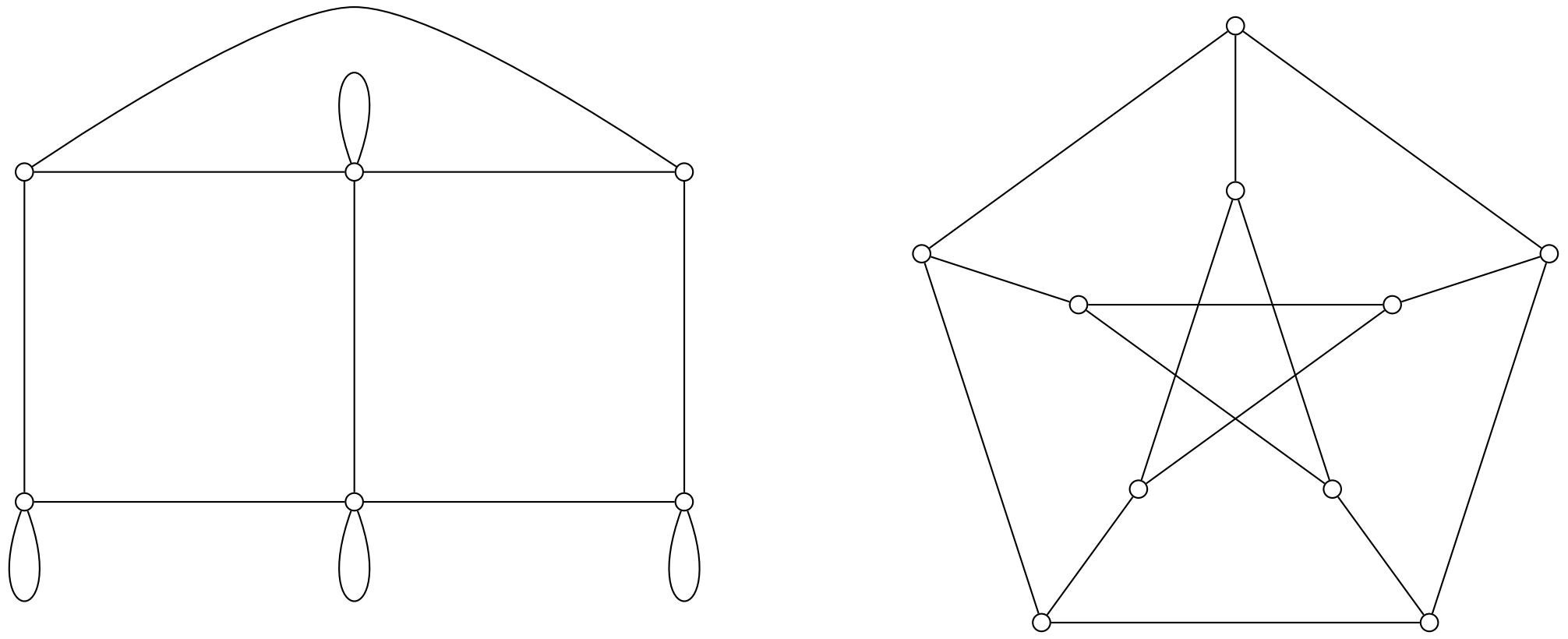

\begin{proof}
The matroid $M(K_1)$ is the empty matroid, realizable by the empty graph. The matroids $M(K_2)$ and $M(K_3)$ are uniform and realizable by Theorem \ref{T.uniform}. The matroids $M(K_4)$ and $M(K_5)$ are realizable via the  graphs depicted in Figure \ref{F.k4k5}.
\end{proof}

\begin{theorem}
\label{T.completegraph}
The graphic matroid $M(K_n)$ is SZF-realizable if and only if $n\in\{1,2,3,4,5\}$.
\end{theorem}

\begin{proof}
This follows immediately from Propositions \ref{P.completeno} and \ref{P.completeyes}.
\end{proof}

\begin{proposition}
\label{P.smallbipartite}
The matroids $M(K_{1,n})$, $M(K_{2,2})$, $M(K_{2,3})$, and $M(K_{3,3})$ are SZF-realizable, for $n\in \mathbb{N}$.
\end{proposition}

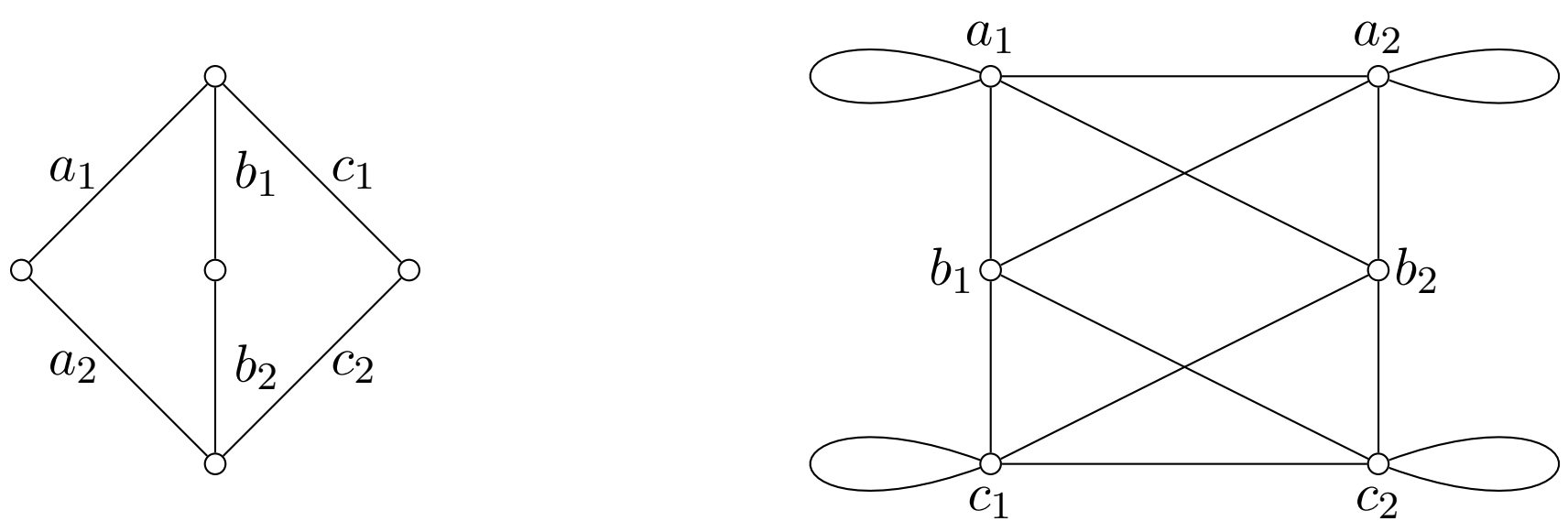
\begin{figure}
    \begin{tikzpicture}[scale=1.4,
  vertex/.style={circle,draw,fill=white,inner sep=1.5pt}
]

%
% K_{2,3} GRAPH
%

% vertices (top, middle, bottom, left, right)
\node[vertex] (t) at (0,2) {};      % top
\node[vertex] (m) at (0,1) {};      % middle
\node[vertex] (b) at (0,0) {};     % bottom
\node[vertex] (l) at (-1,1) {};     % left
\node[vertex] (r) at ( 1,1) {};     % right

% graph label
\node at (0,-0.75) {$K_{2,3}$};

% vertical edges with b_1, b_2
\draw (t) -- node[midway,right] {$b_1$} (m);
\draw (m) -- node[midway,right] {$b_2$} (b);

% left side edges with c_1, c_2
\draw (t) -- node[midway,left] {$a_1$} (l);
\draw (l) -- node[midway,left] {$a_2$} (b);

% right side edges with a_1, a_2
\draw (t) -- node[midway,right] {$c_1$} (r);
\draw (r) -- node[midway,right] {$c_2$} (b);

%
% K_{2,3} SZF GRAPH
%

% Left column
\node[vertex] (L1) at (4,2) {};
\node[vertex] (L2) at (4,1) {};
\node[vertex] (L3) at (4,0) {};

% Right column
\node[vertex] (R1) at (6,2) {};
\node[vertex] (R2) at (6,1) {};
\node[vertex] (R3) at (6,0) {};

% Loops
\draw (L1) to[out=200, in=160, looseness=80] (L1);
\draw (L3) to[out=200, in=160, looseness=80] (L3);
\draw (R1) to[out=340, in=20, looseness=80] (R1);
\draw (R3) to[out=340, in=20, looseness=80] (R3);

% Vertical edges
\draw (L1) -- (L2);
\draw (L2) -- (L3);
\draw (R1) -- (R2);
\draw (R2) -- (R3);

% labels
\node at (4,2.2) {$a_1$};
\node at (3.8,1) {$b_1$};
\node at (4,-.2) {$c_1$};
\node at (6,2.2) {$a_2$};
\node at (6.2,1) {$b_2$};
\node at (6,-.2) {$c_2$};

\node at (5,-0.75) {$M(K_{2,3})_{SZF}$};

% Horizontal edges
\draw (L1) -- (R1);
\draw (L3) -- (R3);

% Diagonal edges
\draw (L1) -- (R2);
\draw (L2) -- (R1);
\draw (L2) -- (R3);
\draw (L3) -- (R2);

\end{tikzpicture}
    \caption{An edge labeling of $K_{2,3}$ and a graph $M(K_{2,3})_{SZF}$ that gives rise to $M(K_{2,3})$ via the SZF closure operator.}
    \label{F.k23}
\end{figure}

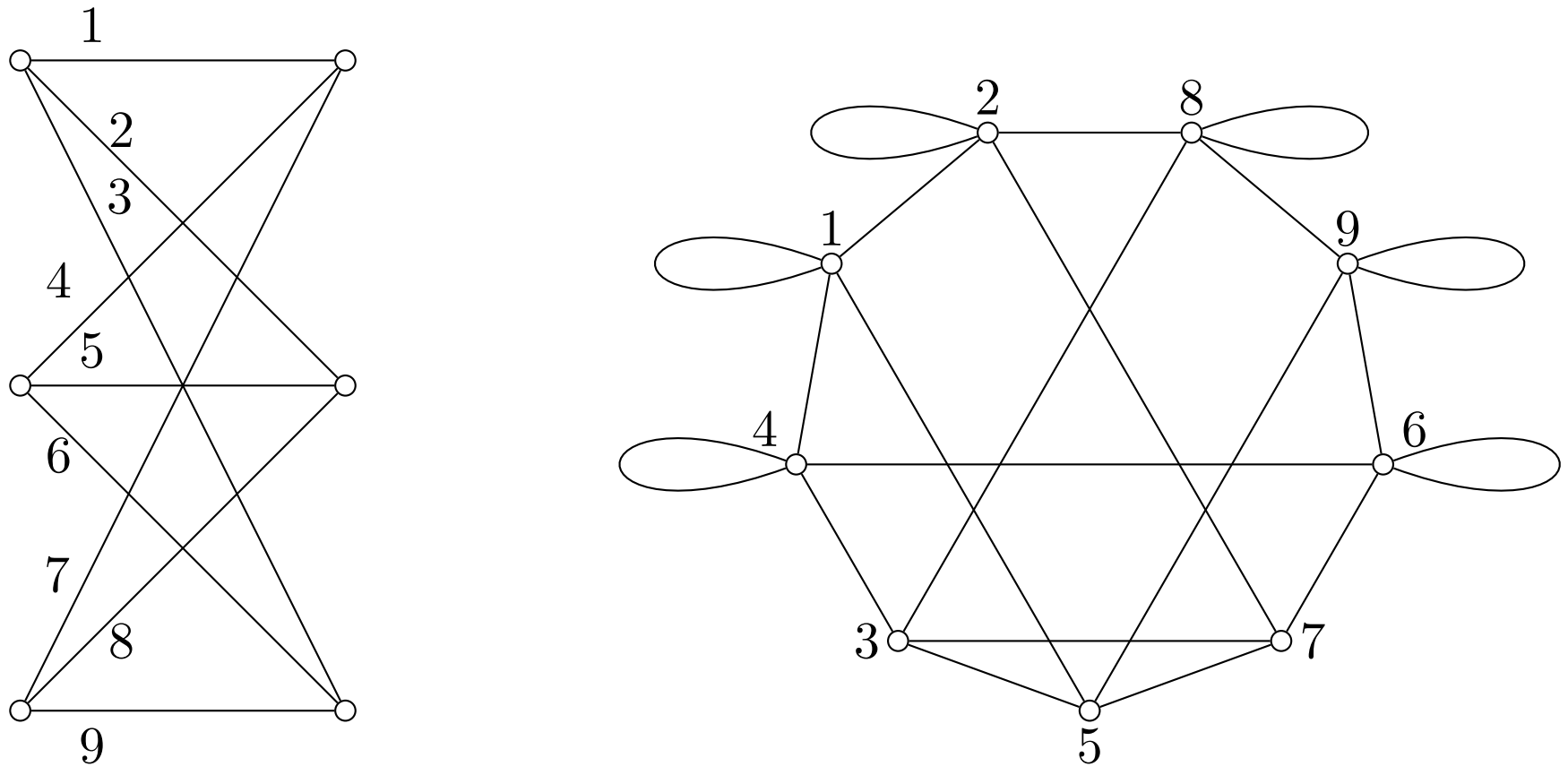
\begin{figure}
      \begin{tikzpicture}[scale=1.2,
  vertex/.style={circle,draw,fill=white,inner sep=1.5pt}
]

% vertices (top, middle, bottom, left, right)
\node[vertex] (lt) at (0,2) {};      % top
\node[vertex] (lm) at (0,0) {};      % middle
\node[vertex] (lb) at (0,-2) {};     % bottom
\node[vertex] (rt) at (2,2) {};     % left
\node[vertex] (rm) at ( 2,0) {};     % right
\node[vertex] (rb) at (2,-2) {};

% edges from top left vertex
\draw (lt) -- node[pos=.2,above] {$1$} (rt);
\draw (lt) -- node[pos=.2, right] {$2$} (rm);
\draw (lt) -- node[pos=.2, right] {$3$} (rb);

% edges from top left vertex
\draw (lm) -- node[pos=.2,above left] {$4$} (rt);
\draw (lm) -- node[pos=.2,above] {$5$} (rm);
\draw (lm) -- node[pos=.2,left] {$6$} (rb);

% edges from top left vertex
\draw (lb) -- node[pos=.2,left] {$7$} (rt);
\draw (lb) -- node[pos=.2, right] {$8$} (rm);
\draw (lb) -- node[pos=.2,below] {$9$} (rb);

\end{tikzpicture} \hspace{0.5in} \begin{tikzpicture}[scale=1.1,
  vertex/.style={circle,draw,fill=white,inner sep=1.5pt}
]

\def\r{2}

% 9 vertices on a circle, with v5 at angle 270 degrees
\foreach \i in {1,...,9} {
  \pgfmathsetmacro{\ang}{270 + 40*(\i-5)}
  \node[vertex] (v\i) at (\ang:\r) {};
}

% Outer 9-cycle
\foreach \i/\j in {1/2,2/3,3/4,4/5,5/6,6/7,7/8,8/9,9/1} {
  \draw (v\i) -- (v\j);
}

% Example loops
\foreach \i in {1,2,3} {
  \draw (v\i) to[out=200, in=160, looseness=80] (v\i);
}

\foreach \i in {7,8,9} {
  \draw (v\i) to[out=340, in=20, looseness=80] (v\i);
}

% Interior edges
\draw (v1) -- (v6);
\draw (v2) -- (v5);
\draw (v4) -- (v9);
\draw (v8) -- (v5);
\draw (v4) -- (v6);
\draw (v3) -- (v7);

% vertex labels
\node[above] at (v1) {2};
\node[above] at (v2) {1};
\node[above left] at (v3) {4};
\node[left] at (v4) {3};
\node[below] at (v5) {5};
\node[right] at (v6) {7};
\node[above right] at (v7) {6};
\node[above] at (v8) {9};
\node[above] at (v9) {8};

\end{tikzpicture}

    \caption{The graph $K_{3,3}$ with edge labels on the left and a graph $M(K_{3,3})_{SZF}$ that gives rise to $M(K_{3,3})$ via the SZF closure operator on the right.}
    \label{F.k33}
\end{figure}

\begin{proof}
    The matroid $M(K_{1,n}) \cong U_{n,n}$ is SZF-realizable by Theorem \ref{T.uniform}. The matroid $M(K_{2,2})$ is SZF-realizable by Proposition \ref{P.cyclematroid}. It is not difficult to check that $M(K_{2,3})_{SZF}$ pictured in Figure \ref{F.k23} gives rise to $M(K_{2,3})$. Similarly, one can check that the graph on the right in Figure \ref{F.k33} gives rise to $M(K_{3,3})$.
\end{proof}

\begin{proposition}
\label{P.2bipartite}
The matroid $M(K_{2,n})$ is not SZF-realizable when $n\geq 4$.
\end{proposition}

\begin{proof}
The graph $K_{2,n}$ has $2n$ edges, each of which is in a series class of size two. From Proposition \ref{P.seriesneighborhoods}, we know that any $M(K_{2,n})_{SZF}$ must have only $n$ distinct neighborhoods (one for each series class), but the matroid has $\binom{n}{2}$ distinct circuit-flats (corresponding to 4-cycles of $K_{2,n}$). Thus, $M(K_{2,n})$ is not SZF-realizable when $n\geq 4$, by Corollary \ref{C.circuitflats}. 
\end{proof}

\begin{proposition}
\label{P.largebipartite}
The matroid $M(K_{m,n})$ is not SZF-realizable when $m,n\geq 3$ and $m+n\geq 7$.
\end{proposition}

\begin{proof}
Observe that $M(K_{m,n})$ has $mn$ elements and $\binom{m}{2} \binom{n}{2}$ circuit-flats. If $m+n\geq 7$ then, by Corollary \ref{C.circuitflats}, the matroid is not SZF-realizable. 
\end{proof}

\begin{theorem}
\label{T.completebipartite}
The graphic matroid $M(K_{m,n})$ is SZF-realizable if and only if
\begin{enumerate}
    \item $m=1$ or $n=1$; or
    \item $m\in \{2,3\}$ and $n\in \{2,3\}$.
\end{enumerate}
\end{theorem}

\begin{proof}
The result follows immediately upon combining Propositions \ref{P.smallbipartite}, \ref{P.2bipartite}, and \ref{P.largebipartite}.    
\end{proof}

For the remainder of this section, the following theorem of Oxley and Whittle \cite{nonspanningCircuits} will be useful, which is restated here for convenience.

\begin{theorem}
\label{T.nonspanningcircuits}
Let $M$ be a matroid on a non-empty set $E$ and $\mathcal{C}'$ be its set of non-spanning circuits. Then $M$ is the unique matroid on $E$ having $\mathcal{C}'$ as its set of non-spanning circuits if and only if the following conditions hold:
\begin{enumerate}
    \item $\mathcal{C}'$ is non-empty and $\max\{|C|: C\in \mathcal{C}'\}=r(M)$.
    \item Either $(i)$ every $(r(M)+1)$-element subset of $E$ contains a member of $\mathcal{C}'$; or $(ii)$ there is an $(r(M)+1)$-element subset $X$ of $E$ and an element $e$ of $E-X$ such that $X$ contains no member of $\mathcal{C}'$ but $X\cup e$ contains at least two members of $\mathcal{C}'$.
\end{enumerate}
\end{theorem}

We have already seen that there are planar graphs whose associated matroids are not SZF-realizable. We now turn our attention to a more restricted class: outerplanar graphs.

\begin{definition}
For a planar graph $G$, the \textit{weak dual} of $G$, denoted $G^{**}$ is obtained from $G$ by taking the dual graph $G^*$ and then deleting the vertex representing the infinite outer face.
\end{definition}

\begin{definition}
    For an embedding of an outerplanar graph $G$, we call the non-infinite faces of $G$ that have at least one edge on the outer cycle of $G$ \textit{boundary faces} and non-infinite faces that have no edges on the outer cycle of $G$ \textit{internal faces}.
\end{definition}

\begin{definition}
    A circuit $C$ in a matroid $M$ is called \textit{chordless} if there exists no circuit $C'$ of $M$ such that $|C'\backslash C|=1$ and $|C'|<|C|$.
\end{definition}

Tutte's \cite{tutte} well-known result that $U_{2,4}$ is the only excluded minor for binary matroids will be used in the following proposition.

\begin{proposition}
\label{P.chordless}
    In a simple binary matroid $M$, the chordless circuits of $M$ are exactly the circuit-flats of $M$.
\end{proposition}

\begin{proof}
    Let $C$ be a circuit-flat of $M$ and assume toward contradiction that it is not chordless. Then there is some circuit $C'$ such that $|C'\backslash C| = 1$ and, by Theorem \ref{T.closure-ccts}, $C$ is not a flat. This contradiction indicates that every circuit-flat of $M$ is a chordless circuit.

    Now assume that $C$ is a chordless circuit of $M$. If $C$ is not a flat then, by Theorem \ref{T.closure-ccts}, there is some circuit $C'$ of $M$ such that $|C'\backslash C|=1$. As $C$ is chordless, we know $|C'|\geq |C|$, which implies $|C'|=|C|$ and $r(C)=r(C')$. As $M$ is a simple binary matroid, this implies that $M$ restricted to $\cl(C)$ has a $U_{n-2,n}$-minor for $n\geq 4$, which contradicts the fact that $M$ is binary.  
\end{proof}

The following theorem from \cite{horns} will be useful in the corollary that follows.

\begin{theorem}
\label{T.horns}
    Let $M$ be a simple binary matroid. Then the set of chordless circuits of $M$ is the unique minimum generator for the set of all circuits of $M$. 
\end{theorem}

\begin{corollary}
\label{C.cctflatouterplanar}
Let $G$ be a simple outerplanar graph. Then the circuit-flats of $M(G)$ contain enough information for an SZF closure operator to completely determine $M(G)$.
\end{corollary}

\begin{proof}
This follows immediately from Proposition \ref{P.chordless}, Theorem \ref{T.horns}, and the fact that all graphic matroids are also binary matroids.
\end{proof}

\begin{corollary}
\label{C.outerplanar}
Let $G$ be a simple outerplanar graph. Then, for any outerplane embedding of $G$, the non-infinite faces of $G$ provide enough information for an SZF closure operator to completely determine $M(G)$.
\end{corollary}

With Corollary \ref{C.outerplanar} in mind, we are ready to prove the main result of the paper, which is restated from the introduction below.

\outerplanar*

\begin{proof}
    Consider an outerplane embedding of $G$ and the associated weak dual $G^{**}$. By Proposition \ref{P.2conn}, we may assume that $G$ is $2$-connected, so $G^{**}$ is a tree and every edge of $G$ appears in the boundary of some face that corresponds to a vertex in $G^{**}$. Consider a maximal path $P_0$ in $G^{**}$ and label its vertices, and the corresponding faces in $G$, with $F_1,\ldots,F_k$. Label the path edges in $P_0$, and their corresponding edges in $G$, with $f_1,\ldots,f_{k-1}$ where $F_i\cap F_{i+1}=\{f_i\}$ for $i\in\{1,\ldots,k-1\}$. For face $F_i$, where $i\in\{2,\ldots,k-1\}$, we call the edges $F_i\setminus\{f_{i-1},f_{i}\}$ the ``leaf edges" of $F_i$ in $P_0$. When $k\geq 2$, we similarly say $F_1\setminus\{f_1\}$ are ``leaf edges" of $F_1$ and $F_k\setminus\{f_{k-1}\}$ are ``leaf edges" of $F_k$. If $k=1$, all edges of $F_1$ are leaf edges.

    We begin constructing $M(G)_{SZF}$, with vertex set $E(G)$, by assigning $n(f_i)=F_{k-i+1}$ for $i\in\{1,\ldots,k-1\}$ and $n(e_k)=F_1$ for all leaf edges $e_k\in F_k$. For $i\in\{1,\ldots, k-1\}$ and leaf edges $e_i\in F_i$, we set $n(e_i)$ to be the set of all edges on the boundary of $F_1\cup F_2 \cup \ldots \cup F_{k-i+1}$. Note that this provides a neighborhood to every edge associated with each $F_i$, maintains a well-defined graph structure, and encodes each $F_i$ as a neighborhood in $M(G)_{SZF}$. Each assigned neighborhood set is a face cycle or a symmetric difference of face cycles and is therefore cyclic in $M(G)$.

    If $G^{**}$ is a path, then all neighborhoods have now been assigned. If $G^{**}$ is not a path, then delete the vertices of $P_0$ and repeat the process above for a path $P_1$ in the remaining forest, maximal with its initial vertex fixed, where the initial vertex of $P_1$ is adjacent in the original $G^{**}$ to a vertex from an earlier path. If an edge $e$ shows up as a leaf edge in multiple paths, then define its neighborhood to be the union of the two sets indicated by the construction above. Continue this process, deleting the vertices of each selected path, until no vertices of remain.

    At each step, the previously selected vertices always induce a connected subtree of $G^{**}$. Then each remaining component has exactly one edge joining it to that subtree. Since each subsequent path starts at one of these attachment vertices and is maximal, it ends at a leaf of the original $G^{**}$. Thus, for each selected path, the terminal face $F_k$ has an outer-boundary edge $e$. Since $e$ belongs to only one face, its assigned neighborhood $n(e)=F_1$ never changes. Similarly, each assigned neighborhood $n(f_i)=F_{k-i+1}$ is never changed since both faces incident with $f_i$ belong to the same path. Thus the edges of each face are assigned to a neighborhood. By Corollary \ref{C.outerplanar}, the neighborhoods generate the closure of $M(G)$. Since each neighborhood is always a union of circuits in $M(G)$ by construction, the resulting set of neighborhoods therefore defines an SZF realization $M(G)_{SZF}$ of $M(G)$.
    \end{proof}

\section*{Declaration of generative AI and AI-assisted technologies in the manuscript preparation process}

During the preparation of this work, the authors used ChatGPT for mathematical exploration and computational code.  Generative AI was not used to write any portion of this manuscript. After using this tool/service, the authors reviewed and edited the content as needed and they take full responsibility for the content of the published article.

\end{document}